\documentclass[psamsfonts]{amsart}

\usepackage{amssymb,amsfonts}
\usepackage[all,arc]{xy}
\usepackage{enumerate}
\usepackage{mathrsfs}
\usepackage{graphicx}

\newtheorem{theorem}{Theorem}[section]
\newtheorem{corollary}[theorem]{Corollary}
\newtheorem{proposition}[theorem]{Proposition}
\newtheorem{lemma}[theorem]{Lemma}

\theoremstyle{definition}
\newtheorem{definition}[theorem]{Definition}

\newtheorem{remark}[theorem]{Remark}

\DeclareMathOperator{\FId}{FId}
\DeclareMathOperator{\Id}{Id}
\DeclareMathOperator{\Max}{Max}

\DeclareMathOperator{\Q}{Q}
\DeclareMathOperator{\MC}{MC}

\DeclareMathOperator{\Frac}{Frac}
\DeclareMathOperator{\F}{F}



\numberwithin{equation}{section}

\begin{document}

\title[\bfseries Invertible Ideals and Gaussian Semirings]{\bfseries Invertible Ideals and Gaussian Semirings}

\author[Shaban Ghalandarzadeh, Peyman Nasehpour, and Rafieh Razavi]{Shaban Ghalandarzadeh, Peyman Nasehpour, and Rafieh Razavi}

\address{Shaban Ghalandarzadeh, Faculty of Mathematics, K. N. Toosi University of Technology, Tehran, Iran}
\email{ghalandarzadeh@kntu.ac.ir}

\address{Peyman Nasehpour, Department of Engineering Science, Golpayegan University of Technology, Golpayegan, Iran}
\address{Peyman Nasehpour, Department of Engineering Science, Faculty of Engineering, University of Tehran, Tehran, Iran}

\email{nasehpour@gmail.com, nasehpour@gut.ac.ir}

\address{Rafieh Razavi, Faculty of Mathematics, K. N. Toosi University of Technology, Tehran, Iran}
\email{rrazavi@mail.kntu.ac.ir}



\begin{abstract}
In the first section, we introduce the notions of fractional and invertible ideals of semirings and characterize invertible ideals of a semidomain. In section two, we define Pr\"{u}fer semirings and characterize them in terms of valuation semirings. In this section, we also characterize Pr\"{u}fer semirings in terms of some identities over its ideals such as $(I + J)(I \cap J) = IJ$ for all ideals $I$, $J$ of $S$. In the third section, we give a semiring version for the Gilmer-Tsang Theorem, which states that for a suitable family of semirings, the concepts of Pr\"{u}fer and Gaussian semirings are equivalent. At last, we end this paper by giving a plenty of examples for proper Gaussian and Pr\"{u}fer semirings.
\end{abstract}

\maketitle

\section{Introduction}

Vandiver introduced the term ``semi-ring" and its structure in 1934 \cite{Vandiver1934}, though the early examples of semirings had appeared in the works of Dedekind in 1894, when he had been working on the algebra of the ideals of commutative rings \cite{Dedekind1894}. Despite the great efforts of some mathematicians on semiring theory in 1940s, 1950s, and early 1960s, they were apparently not successful to draw the attention of mathematical society to consider the semiring theory as a serious line of mathematical research. Actually, it was in the late 1960s that semiring theory was considered a more important topic for research when real applications were found for semirings. Eilenberg and a couple of other mathematicians started developing formal languages and automata theory systematically \cite{Eilenberg1974}, which have strong connections to semirings. Since then, because of the wonderful applications of semirings in engineering, many mathematicians and computer scientists have broadened the theory of semirings and related structures \cite{Glazek2002} and \cite{HebischWeinert1998}. As stated in \cite[p. 6]{Golan1999(b)}, multiplicative ideal theoretic methods in ring theory are certainly one of the major sources of inspiration and problems for semiring theory. In the present paper, we develop some ring theoretic methods of multiplicative ideal theory for semirings as follows:

Let, for the moment, $R$ be a commutative ring with a nonzero identity. The Dedekind-Mertens lemma in ring theory states that if $f$ and $g$ are two elements of the polynomial ring $R[X]$, then there exists a natural number $n$ such that $ c(f)^{n-1} c(fg) = c(f)^n c(g)$, where by the content $c(f)$ of an arbitrary polynomial $ f \in R[X]$, it is meant the $R$-ideal generated by the coefficients of $f$. From this, it is clear that if $R$ is a Pr\"{u}fer domain, then $R$ is Gaussian, i.e. $c(fg) = c(f)c(g)$ for all $f,g \in R[X]$.

Gilmer in \cite{Gilmer1967} and Tsang in \cite{Tsang1965}, independently, proved that the inverse of the above statement is also correct in this sense that if $R$ is a Gaussian domain, then $R$ is a Pr\"{u}fer domain.

Since Gaussian semirings were introduced in Definition 7 in \cite{Nasehpour2016} and the Dedekind-Mertens lemma was proved for subtractive semirings in Theorem 3 in the same paper, our motivation for this work was to see how one could define invertible ideals for semirings to use them in Dedekind-Mertens lemma and discover another family of Gaussian semirings. We do emphasize that the definition of Gaussian semiring used in our paper is different from the one investigated in \cite{DalePitts1978} and \cite{HebischWeinert1987}. We also asked ourselves if some kind of a Gilmer-Tsang Theorem held for polynomial semirings. Therefore, we were not surprised to see while investigating these questions, we needed to borrow some definitions and techniques - for example Pr\"{u}fer domains and a couple of other concepts mentioned in \cite{LarsenMcCarthy1971} and \cite{Gilmer1972} - from multiplicative ideal theory for rings . In most cases, we also constructed examples of proper semirings - semirings that are not rings - satisfying the conditions of those definitions and results to show that what we bring in this paper are really generalizations of their ring version ones. Since different authors have used the term ``semiring" with some different meanings, it is essential, from the beginning, to clarify what we mean by a semiring.

In this paper, by a semiring, we understand an algebraic structure, consisting of a nonempty set $S$ with two operations of addition and multiplication such that the following conditions are satisfied:

\begin{enumerate}
\item $(S,+)$ is a commutative monoid with identity element $0$;
\item $(S,\cdot)$ is a commutative monoid with identity element $1 \not= 0$;
\item Multiplication distributes over addition, i.e. $a(b+c) = ab + ac$ for all $a,b,c \in S$;
\item The element $0$ is the absorbing element of the multiplication, i.e. $s\cdot 0=0$ for all $s\in S$.
\end{enumerate}

From the above definition, it is clear for the reader that the semirings are fairly interesting generalizations of the two important and broadly studied algebraic structures, i.e. rings and bounded distributive lattices.

A nonempty subset $I$ of a semiring $S$ is defined to be an ideal of $S$ if $a, b \in I$ and $ s \in S$ implies that $a + b, sa \in I $ \cite{Bourne1951}. An ideal $I$ of a semiring $S$ is said to be subtractive, if $a+b \in I$ and $a \in I$ implies $b \in I$ for all $a,b \in S$. A semiring $S$ is said to be subtractive if every ideal of the semiring $S$ is subtractive. An ideal $P$ of $S$ is called a prime ideal of $S$ if $P \neq S$ and $ab\in P$ implies that $a\in P$ or $b\in P$ for all $a,b \in S$.\\

In $\S1$, we define fractional and invertible ideals and show that any invertible ideal of a local semidomain is principal (See Definitions \ref{fractional} and \ref{invertible} and Proposition \ref{invertibleisprincipal1}). Note that a semiring $S$ is called a semidomain if for any nonzero element $s$ of $S$, $sb = sc$ implies that $b=c$. A semiring is said to be local if it has only one maximal ideal. We also prove that any invertible ideal of a weak Gaussian semi-local semidomain is principal (See Theorem \ref{invertibleisprincipal2}). Note that a semiring is defined to be a weak Gaussian semiring if each prime ideal of the semiring is subtractive \cite[Definition 18]{Nasehpour2016} and a semiring is said to be semi-local if the set of its maximal ideals is finite. Also, note that localization of semirings has been introduced and investigated in \cite{Kim1985}. It is good to mention that an equivalent definition for the localization of semirings has been given in \cite[\S11]{Golan1999(b)}.

At last, in Theorem \ref{invertiblelocalization2}, we show that if $I$ is a nonzero finitely generated ideal of a semidomain $S$, then $I$ is invertible if and only if $I_{\mathfrak{m}}$ is a principal ideal of $S_{\mathfrak{m}}$ for each maximal ideal $\mathfrak{m}$ of $S$.\\

In $\S2$, we observe that if $S$ is a semiring, then every nonzero finitely generated ideal of $S$ is an invertible ideal of $S$ if and only if every nonzero principal and every nonzero 2-generated ideal of $S$ is an invertible ideal of $S$ (Check Theorem \ref{prufer}). This result and a nice example of a proper semiring having this property, motivate us to define Pr\"{u}fer semiring, the semiring that each of its nonzero finitely generated ideals is invertible (See Definition \ref{prufer-def}). After that, in Theorem \ref{prufer2}, we prove that a semidomain $S$ is a Pr\"{u}fer semiring if and only if one of the following equivalent statements holds:

\begin{enumerate}

\item $I(J\cap K)=IJ\cap IK$ for all ideals $I$, $J$, and $K$ of $S$,
\item $(I+J)(I\cap J)=IJ$ for all ideals $I$ and $J$ of $S$,
\item $[(I+J):K]=[I:K]+[J:K]$ for all ideals $I$, $J$, and $K$ of $S$ with $K$ finitely generated,
\item $[I:J]+[J:I]=S$ for all finitely generated ideals $I$ and $J$ of $S$,
\item $[K:I\cap J]=[K:I]+[K:J]$ for all ideals $I$, $J$, and $K$ of $S$ with $I$ and $J$ finitely generated.

\end{enumerate}

Note that, in the above, it is defined that $[I \colon J]=\lbrace s\in S \colon sJ \subseteq I \rbrace$. Also, note that this theorem is the semiring version of Theorem 6.6 in \cite{LarsenMcCarthy1971}, though we give partly an alternative proof for the semiring generalization of its ring version.

In $\S2$, we also characterize Pr\"{u}fer semirings in terms of valuation semirings. Let us recall that a semidomain is valuation if its ideals are totally ordered by inclusion \cite[Theorem 2.4]{Nasehpour2018}. In fact, in Theorem \ref{valuationlocalization}, we prove that a semiring $S$ is Pr\"{u}fer if and only if one of the following statements holds:

\begin{enumerate}

\item For any prime ideal $\mathfrak{p}$ of $S$, $S_\mathfrak{p}$ is a valuation semidomain.

\item For any maximal ideal $\mathfrak{m}$ of $S$, $S_\mathfrak{m}$ is a valuation semidomain.

\end{enumerate}

A nonzero ideal $I$ of a semiring $S$ is called a cancellation ideal, if $IJ=IK$ implies $J=K$ for all ideals $J$ and $K$ of $S$ \cite{LaGrassa1995}. Let $f\in S[X]$ be a polynomial over the semiring $S$. The content of $f$, denoted by $c(f)$, is defined to be the $S$-ideal generated by the coefficients of $f$. It is, then, easy to see that $c(fg) \subseteq c(f)c(g)$ for all $f,g \in S[X]$. Finally, a semiring $S$ is defined to be a Gaussian semiring if $c(fg)=c(f)c(g)$ for all $f,g \in S[X]$ \cite[Definition 8]{Nasehpour2016}.\\

In $\S3$, we discuss Gaussian semirings and prove a semiring version of the Gilmer-Tsang Theorem with the following statement (See Theorem \ref{Gilmer-Tsang}):

Let $S$ be a subtractive semiring such that every nonzero principal ideal of $S$ is invertible and $ab\in (a^2,b^2)$ for all $a,b \in S$. Then the following statements are equivalent:

 \begin{enumerate}

 \item $S$ is a Pr\"{u}fer semiring,
 \item Each nonzero finitely generated ideal of $S$ is cancellation,
 \item $[IJ \colon I]=J$ for all ideals $I$ and $J$ of $S$,
 \item $S$ is a Gaussian semiring.
 \end{enumerate}

At last, we end this paper by giving a plenty of examples of proper Gaussian and Pr\"{u}fer semirings in Theorem \ref{notring} and Corollary \ref{notring2}. Actually, we prove that if $S$ is a Pr\"{u}fer semiring (say for example $S$ is a Pr\"{u}fer domain), then $\FId(S)$ is a Pr\"{u}fer semiring, where by $\FId(S)$ we mean the semiring of finitely generated ideals of $S$.

In this paper, all semirings are assumed to be commutative with a nonzero identity. Unless otherwise stated, our terminology and notation will follow as closely as possible that of \cite{Gilmer1972}.

\section{Fractional and invertible ideals of semirings}

In this section, we introduce fractional and invertible ideals for semirings and prove a couple of interesting results for them. Note that whenever we feel it is necessary, we recall concepts related to semiring theory to make the paper as self-contained as possible.

Let us recall that a nonempty subset $I$ of a semiring $S$ is defined to be an ideal of $S$ if $a, b \in I$ and $ s \in S$ implies that $a + b, sa \in I $ \cite{Bourne1951}. Also, $T \subseteq S$ is said to be a multiplicatively closed set of $S$ provided that if $a,b \in T$, then $ab\in T$. The localization of $S$ at $T$ is defined in the following way:

First define the equivalent relation $ \sim $ on $S \times T$ by $(a,b) \sim (c,d)$, if $tad=tbc$ for some $t\in T$. Then Put $S_T$ the set of all equivalence classes of $S \times T$ and define addition and multiplication on $S_T$ respectively by $[a,b]+[c,d] = [ad+bc, bd]$ and $[a,b]\cdot[c,d] = [ac,bd]$, where by $[a,b]$, also denoted by $a/b$, we mean the equivalence class of $(a,b)$. It is, then, easy to see that $S_T$ with the mentioned operations of addition and multiplication in above is a semiring \cite{Kim1985}.

Also, note that an element $s$ of a semiring $S$ is said to be multiplicatively-cancellable (abbreviated as MC), if $sb=sc$ implies $b=c$ for all $b,c \in S$. For more on MC elements of a semiring, refer to \cite{ElBashirHurtJancarikKepka2001}. We denote the set of all MC elements of $S$ by $\MC(S)$. It is clear that $\MC(S)$ is a multiplicatively closed set of $S$. Similar to ring theory, total quotient semiring $\Q(S)$ of the semiring $S$ is defined as the localization of $S$ at $\MC(S)$. Note that $\Q(S)$ is also an $S$-semimodule. For a definition and a general discussion of semimodules, refer to \cite[\S 14]{Golan1999(b)}. Now, we define fractional ideals of a semiring as follows:

\begin{definition}

\label{fractional}

\textbf{Fractional ideal}. We define a fractional ideal of a semiring $S$ to be a subset $I$ of the total quotient semiring $\Q(S)$ of $S$ such that:

\begin{enumerate}

\item I is an $S$-subsemimodule of $\Q(S)$, that is, if $a, b \in I$ and $ s \in S$, then $a + b \in I $ and $sa \in I$.

\item There exists an MC element $d\in S$ such that $dI \subseteq S$.

\end{enumerate}

\end{definition}

Let us denote the set of all nonzero fractional ideals of $S$ by $\Frac(S)$. It is easy to check that $\Frac(S)$ equipped with the following multiplication of fractional ideals is a commutative monoid:

$$I\cdot J = \{a_1 b_1 + \cdots + a_n b_n : a_i \in I, b_i \in J \}.$$

\begin{definition}

\label{invertible}

\textbf{Invertible ideal}. We define a fractional ideal $I$ of a semiring $S$ to be invertible if there exists a fractional ideal $J$ of $S$ such that $IJ=S$.

\end{definition}

Note that if a fractional ideal $I$ of a semiring $S$ is invertible and $IJ =S$, for some fractional ideal $J$ of $S$, then $J$ is unique and we denote that by $I^{-1}$. It is clear that the set of invertible ideals of a semiring equipped with the multiplication of fractional ideals is an Abelian group.

\begin{theorem}

\label{frac2}

Let $S$ be a semiring with its total quotient semiring $\Q(S)$.

\begin{enumerate}

\item If $I \in \Frac(S)$ is invertible, then $I$ is a finitely generated $S$-subsemimodule of $\Q(S)$.

\item If $I, J \in \Frac(S) $ and $ I \subseteq J$ and $J$ is invertible, then there is an ideal $K$ of $S$ such that $I = JK$.

\item If $I \in \Frac(S)$, then $I$ is invertible if and only if there is a fractional ideal $J$ of $S$ such that $IJ$ is principal and generated by an MC element of $\Q(S)$.

\end{enumerate}

\begin{proof}
The proof of this theorem is nothing but the mimic of the proof of its ring version in \cite[Proposition 6.3]{LarsenMcCarthy1971}.
\end{proof}

\end{theorem}

Let us recall that a semiring $S$ is defined to be a semidomain, if each nonzero element of the semiring $S$ is an MC element of $S$.

\begin{proposition}
Let $S$ be a semiring and $a\in S$. Then the following statements hold:

\begin{enumerate}

\item The principal ideal $(a)$ is invertible if and only if $a$ is an MC element of $S$.

\item The semiring $S$ is a semidomain if and only if each nonzero principal ideal of $S$ is an invertible ideal of $S$.

\end{enumerate}

\begin{proof}
Straightforward.
\end{proof}

\end{proposition}

Prime and maximal ideals of a semiring are defined similar to rings (\cite[\S7]{Golan1999(b)}). Note that the set of the unit elements of a semiring $S$ is denoted by $U(S)$. Also note that when $S$ is a semidomain, $\MC(S) = S-\{0\}$ and the localization of $S$ at $\MC(S)$ is called the semifield of fractions of the semidomain $S$ and usually denoted by $\F(S)$ \cite[p. 22]{Golan1999(a)}.

\begin{proposition}

\label{invertibleisprincipal1}

Any invertible ideal of a local semidomain is principal.

\begin{proof}

Let $I$ be an invertible ideal of a local semidomain $(S, \mathfrak{m})$. It is clear that there are $s_1,\ldots,s_n \in S$ and $t_1,\ldots,t_n \in \F(S)$, such that $s_1 t_1 + \cdots + s_n t_n = 1$. This implies that at least one of the elements $s_i t_i$ is a unit, since if all of them are nonunit, their sum will be in $\mathfrak{m}$ and cannot be equal to $1$. Assume that $s_1 t_1 \in U(S)$. Now we have $S = (s_1)(t_1) \subseteq I (t_1) \subseteq I I^{-1} = S$, which obviously implies that $I = (s_1)$ and the proof is complete.
\end{proof}

\end{proposition}

Let us recall that an ideal $I$ of a semiring $S$ is said to be subtractive, if $a+b \in I$ and $a \in I$ implies $b \in I$ for all $a,b \in S$. Now we prove a similar statement for weak Gaussian semirings introduced in \cite{Nasehpour2016}. Note that any prime ideal of a weak Gaussian semiring is subtractive (\cite[Theorem 19]{Nasehpour2016}). Using this property, we prove the following theorem:

\begin{theorem}

\label{invertibleisprincipal2}

Any invertible ideal of a weak Gaussian semi-local semidomain is principal.

\begin{proof}
Let $S$ be a weak Gaussian semi-local semidomain and $\Max(S) = \{\mathfrak{m}_1, \ldots, \mathfrak{m}_n\}$ and $II^{-1} =S$. Similar to the proof of Proposition \ref{invertibleisprincipal1}, for each $1 \leq i \leq n$, there exist $a_i \in I$ and $b_i \in I^{-1}$ such that $a_i b_i \notin \mathfrak{m}_i$. Since by \cite[Corollary 7.13]{Golan1999(b)} any maximal ideal of a semiring is prime, one can easily check that any $\mathfrak{m}_i$ cannot contain the intersection of the remaining maximal ideals of $S$. So for any $1 \leq i \leq n$, one can find some $u_i$, where $u_i$ is not in $\mathfrak{m}_i$, while it is in all the other maximal ideals of $S$. Put $v = u_1 b_1 + \cdots + u_n b_n$. It is obvious that $v\in I^{-1}$, which causes $vI$ to be an ideal of $S$. Our claim is that $vI$ is not a subset of any maximal ideal of $S$. In contrary assume that $vI$ is a subset of a maximal ideal, say $\mathfrak{m}_1$. This implies that $v a_1\in  \mathfrak{m}_1$. But $$v a_1 = (u_1 b_1 + \cdots + u_n b_n)a_1 .$$

Also note that $u_i b_i a_i \in \mathfrak{m}_1$ for any $i \geq2$. Since $\mathfrak{m}_1$ is subtractive, $u_1 b_1 a_1 \in \mathfrak{m}_1$, a contradiction. From all we said we have that $vI = S$ and finally $I = (v^{-1})$, as required.
\end{proof}

\end{theorem}

The proof of the following lemma is straightforward, but we bring it only for the sake of reference.

\begin{lemma}

\label{invertiblelocalization1}

Let $I$ be an invertible ideal in a semidomain $S$ and $T$ a multiplicatively closed set. Then $I_T$ is an invertible ideal of $S_T$.

\begin{proof}
Straightforward.
\end{proof}

\end{lemma}

Let us recall that if $\mathfrak{m}$ is a maximal ideal of $S$, then $S-\mathfrak{m}$ is a multiplicatively closed set of $S$ and the localization of $S$ at $S-\mathfrak{m}$ is simply denoted by $S_{\mathfrak{m}}$ \cite{Kim1985}. Now, we prove the following theorem:

\begin{theorem}

\label{invertiblelocalization2}

Let $I$ be a nonzero finitely generated ideal of a semidomain $S$. Then $I$ is invertible if and only if $I_{\mathfrak{m}}$ is a principal ideal of $S_{\mathfrak{m}}$ for each maximal ideal $\mathfrak{m}$ of $S$.

\begin{proof} Let $S$ be a semidomain and $I$ a nonzero finitely generated ideal of $S$.

$(\rightarrow):$ If $I$ is invertible, then by Lemma \ref{invertiblelocalization1}, $I_{\mathfrak{m}}$ is invertible and therefore, by Proposition \ref{invertibleisprincipal1}, is principal.

$(\leftarrow):$ Assume that $I_{\mathfrak{m}}$ is a principal ideal of $S_{\mathfrak{m}}$ for each maximal ideal $\mathfrak{m}$ of $S$. For the ideal $I$, define $J := \{x\in \F(S) : xI \subseteq S\}$. It is easy to check that $J$ is a fractional ideal of $S$ and $IJ \subseteq S$ is an ideal of $S$. Our claim is that $IJ = S$. On the contrary, suppose that $IJ \neq S$. So $IJ$ lies under a maximal ideal $ \mathfrak{m}$ of $S$. By hypothesis $I_{\mathfrak{m}}$ is principal. We can choose a generator for $I_{\mathfrak{m}}$ to be an element $z\in I$. Now let $a_1,\ldots, a_n$ be generators of $I$ in $S$. It is, then, clear that for any $a_i$, one can find an $s_i \in S- \mathfrak{m}$ such that $a_i s_i \in (z)$. Set $s = s_1 \cdots s_n$. Since $(sz^{-1}) a_i \in S$, by definition of $J$, we have $sz^{-1} \in J$. But now $s = (sz^{-1}) z \in \mathfrak{m}$, contradicting that $s_i \in S- \mathfrak{m}$ and the proof is complete.
\end{proof}

\end{theorem}

Now the question arises if there is any proper semiring, which each of its nonzero finitely generated ideals is invertible. The answer is affirmative and next section is devoted to such semirings.

\section{Pr\"{u}fer semirings}

The purpose of this section is to introduce the concept of Pr\"{u}fer semirings and investigate some of their properties. We start by proving the following important theorem, which in its ring version can be found in \cite[Theorem 6.6]{LarsenMcCarthy1971}.

\begin{theorem}

\label{prufer}

Let $S$ be a semiring. Then the following statements are equivalent:

\begin{enumerate}

\item Each nonzero finitely generated ideal of $S$ is an invertible ideal of $S$,

\item The semiring $S$ is a semidomain and every nonzero 2-generated ideal of $S$ is an invertible ideal of $S$.

\end{enumerate}

\begin{proof}
Obviously the first assertion implies the second one. We prove that the second assertion implies the first one. The proof is by induction. Let $n>2$ be a natural number and suppose that all nonzero ideals of $S$ generated by less than $n$ generators are invertible ideals and $L=(a_1,a_2,\ldots,a_{n-1},a_n)$ be an ideal of $S$. If we put $I=(a_1)$, $J=(a_2,\ldots,a_{n-1})$ and $K=(a_n)$, then by induction's hypothesis the ideals $I+J$, $J+K$ and $K+I$ are all invertible ideals. On the other hand, a simple calculation shows that the identity $(I+J)(J+K)(K+I)=(I+J+K)(IJ+JK+KI)$ holds. Also since product of fractional ideals of $S$ is invertible if and only if every factor of this product is invertible, the ideal $I+J+K=L$ is invertible and the proof is complete.
\end{proof}

\end{theorem}

A ring $R$ is said to be a Pr\"{u}fer domain if every nonzero finitely generated ideal of $R$ is invertible. It is, now, natural to ask if there is any proper semiring $S$ with this property that every nonzero finitely generated ideal of $S$ is invertible. In the following remark, we give such an example.

\begin{remark}

\label{pruferexample}

Example of a proper semiring with this property that every nonzero finitely generated ideal of $S$ is invertible: Obviously $(\Id(\mathbb Z),+,\cdot)$ is a semidomain, since any element of $\Id(\mathbb Z)$ is of the form $(n)$ such that $n$ is a nonnegative integer and $(a)(b)=(ab)$, for any $a,b \geq 0$. Let $I$ be an arbitrary ideal of $\Id(\mathbb Z)$. Define $A_I$ to be the set of all positive integers $n$ such that $(n) \in I$ and put $m = \min A_I$. Our claim is that $I$ is the principal ideal of $\Id(\mathbb Z)$, generated by $(m)$, i.e. $I=((m))$. For doing so, let $(d)$ be an element of $I$. But then $(\gcd(d,m))= (d)+(m)$ and therefore, $(\gcd(d,m)) \in I$. This means that $ m \leq \gcd(d,m)$, since $m = \min A_I$, while $\gcd(d,m) \leq m$ and this implies that $\gcd(d,m) = m$ and so $m$ divides $d$ and therefore, there exists a natural number $r$ such that $d=rm$. Hence, $(d)=(r)(m)$ and the proof of our claim is finished. From all we said we learn that each ideal of the semiring $\Id(\mathbb Z)$ is a principal and, therefore, an invertible ideal, while obviously it is not a ring.

\end{remark}

By Theorem \ref{prufer} and the example given in Remark \ref{pruferexample}, we are inspired to give the following definition:

\begin{definition}

\label{prufer-def}

We define a semiring $S$ to be a Pr\"{u}fer semiring if every nonzero finitely generated ideal of $S$ is invertible.

\end{definition}

First we prove the following interesting results:

\begin{lemma}

\label{prufer1r}

Let $S$ be a Pr\"{u}fer semiring. Then $I\cap(J+K)=I\cap J+I\cap K$ for all ideals $I$, $J$, and $K$ of $S$.

\begin{proof}

Let $s\in I\cap(J+K)$. So there are $s_1\in J$ and $s_2\in K$ such that $s=s_1+s_2\in I$. If we put $L=(s_1,s_2)$, by definition, we have $LL^{-1}=S$. Consequently, there are $t_1,t_2\in L^{-1}$ such that $s_1t_1+s_2t_2=1$. So $s=ss_1t_1+ss_2t_2$. But $st_1,st_2\in S$, since $s=s_1+s_2\in L$. Therefore, $ss_1t_1\in J$ and $ss_2t_2\in K$. Moreover $s_1t_1,s_2t_2 \in S$ and therefore, $ss_1t_1 ,ss_2t_2\in I$. This implies that $ss_1t_1\in I\cap J$, $ss_2t_2\in I\cap K$, and $s\in I\cap J+I\cap K$, which means that $I\cap(J+K)\subseteq I\cap J+I\cap K$. Since the reverse inclusion is always true, $I\cap(J+K)=I\cap J+I\cap K$ and this finishes the proof.
\end{proof}
\end{lemma}

\begin{lemma}

\label{can1}

Let $S$ be a Pr\"{u}fer semiring. Then the following statements hold:

\begin{enumerate}

\item If $I$ and $K$ are ideals of $S$, with $K$ finitely generated, and if $ I \subseteq K$, then there is an ideal $J$ of $S$ such that $I=JK$.

\item If $IJ=IK$, where $I$, $J$ and $K$ are ideals of $S$ and $I$ is finitely generated and nonzero, then $J=K$.

\end{enumerate}

\begin{proof}
By considering Theorem \ref{frac2}, the assertion $(1)$ holds. The assertion $(2)$ is straightforward.
\end{proof}

\end{lemma}

Note that the second property in Lemma \ref{can1} is the concept of cancellation ideal for semirings, introduced in \cite{LaGrassa1995}:

\begin{definition}

\label{cancelleationidealsdef}

A nonzero ideal $I$ of a semiring $S$ is called a cancellation ideal, if $IJ=IK$ implies $J=K$ for all ideals $J$ and $K$ of $S$.

\end{definition}

\begin{remark}

It is clear that each invertible ideal of a semiring is cancellation. Also, each finitely generated nonzero ideal of a Pr\"{u}fer semiring is cancellation. For a general discussion on cancellation ideals in rings, refer to \cite{Gilmer1972} and for generalizations of this concept in module and ring theory, refer to \cite{NaoumMijbass1997} and \cite{NassehpourYassemi2000}.

\end{remark}

While the topic of cancellation ideals is interesting by itself, we do not go through them deeply. In fact in this section, we only prove the following result for cancellation ideals of semirings, since we need it in the proof of Theorem \ref{Gilmer-Tsang}. Note that similar to ring theory, for any ideals $I$ and $J$ of a semiring $S$, it is defined that $$[I \colon J]=\lbrace s\in S \colon sJ \subseteq I \rbrace.$$ Also, we point out that this result is the semiring version of an assertion mentioned in \cite[Exercise. 4, p. 66]{Gilmer1972}:

\begin{proposition}

\label{can2}

Let $S$ be a semiring and $I$ be a nonzero ideal of $S$. Then the following statements are equivalent:

\begin{enumerate}

\item $I$ is a cancellation ideal of $S$,

\item $[IJ \colon I]=J$ for any ideal $J$ of $S$,

\item $IJ \subseteq IK$ implies $J \subseteq K$ for all ideals $J,K$ of $S$.

\end{enumerate}

\begin{proof}
By considering this point that the equality $[IJ \colon I]I=IJ$ holds for all ideals $I,J$ of $S$, it is then easy to see that (1) implies (2). The rest of the proof is straightforward.
\end{proof}

\end{proposition}

Now we prove an important theorem that is rather the semiring version of Theorem 6.6 in \cite{LarsenMcCarthy1971}. While some parts of our proof is similar to those ones in Theorem 6.6 in \cite{LarsenMcCarthy1971} and Proposition 4 in \cite{S}, other parts of the proof are apparently original.

\begin{theorem}

\label{prufer2}

Let $S$ be a semidomain. Then the following statements are equivalent:

\begin{enumerate}

\item The semiring $S$ is a Pr\"{u}fer semiring,
\item $I(J\cap K)=IJ\cap IK$ for all ideals $I$, $J$, and $K$ of $S$,
\item $(I+J)(I\cap J)=IJ$ for all ideals $I$ and $J$ of $S$,
\item $[(I+J):K]=[I:K]+[J:K]$ for all ideals $I$, $J$, and $K$ of $S$ with $K$ finitely generated,
\item $[I:J]+[J:I]=S$ for all finitely generated ideals $I$ and $J$ of $S$,
\item $[K:I\cap J]=[K:I]+[K:J]$ for all ideals $I$, $J$, and $K$ of $S$ with $I$ and $J$ finitely generated.

\end{enumerate}

\begin{proof}

$(1) \rightarrow (2)$: It is clear that $I(J\cap K)\subseteq IJ\cap IK$. Let $s\in IJ\cap IK$. So we can write $s=\sum_{i=1}^{m}t_iz_i=\sum_{j=1}^{n}t'_jz'_j$ , where $t_i, t'_j\in I$, $z_i \in J$, and $z'_j \in K$ for all $1\leq i \leq m$ and $1\leq j \leq n$. Put $I_1=(t_1,\ldots,t_m)$, $I_2=(t'_1,\ldots,t'_n)$, $J'=(z_1,\ldots,z_m)$, $K'=(z'_1,\ldots,z'_n)$, and $I_3=I_1+I_2$. Then $I_1J'\cap I_2K'\subseteq I_3J'\cap I_3K'\subseteq I_3$. Since $I_3$ is a finitely generated ideal of $S$, by Lemma \ref{can1}, there exists an ideal $L$ of $S$ such that $I_3J'\cap I_3K'=I_3L$. Note that $L={I_3}^{-1}(I_3J'\cap I_3K')\subseteq {I_3}^{-1}(I_3J')=J'$. Moreover $L={I_3}^{-1}(I_3J'\cap I_3K')\subseteq {I_3}^{-1}(I_3K')=K'$. Therefore, $L\subseteq J'\cap K'$. Thus $s\in  I_3J'\cap I_3K'={I_3} L\subseteq I_3(J'\cap K')\subseteq I(J\cap K)$.

$(2) \rightarrow (3)$: Let $I,J\subseteq S$. Then $(I+J)(I\cap J)=(I+J)I\cap(I+J)J \supseteq IJ$. Since the reverse inclusion always holds, $(I+J)(I\cap J)=IJ$.

$(3) \rightarrow (1)$: By hypothesis, every two generated ideal $I = (s_1, s_2)$ is a factor of the invertible ideal $(s_1 s_2)$ and therefore, it is itself invertible. Now by considering Theorem \ref{prufer}, it is clear that the semiring $S$ is a Pr\"{u}fer semiring.

$(1) \rightarrow (4)$: Let $s\in S$ such that $sK\subseteq I+J$. So $sK\subseteq (I+J)\cap K$. By Lemma \ref{prufer1r}, $sK\subseteq I\cap K+J\cap K$. Therefore, $s\in (I\cap K)K^{-1}+(J\cap K)K^{-1}$. Thus $s=\sum_{i=1}^{m}t_iz_i+\sum_{j=1}^{n}t'_jz'_j$, where $z_i, z'_j\in K^{-1}$, $t_i \in I\cap K$, and $t'_j \in J\cap K$ for all $1\leq i \leq m$ and $1\leq j \leq n$. Let $x\in K$ and $1\leq i \leq m$. Then $z_ix,z_it_i\in S$ and so $t_iz_ix\in I\cap K$. Therefore, $(\sum_{i=1}^{m}t_iz_i)K\subseteq I\cap K\subseteq I$. In a similar way, $(\sum_{j=1}^{n}t'_jz'_j)K\subseteq J\cap K\subseteq J$. Thus $s\in [I:K]+[J:K]$. Therefore, $[(I+J):K]\subseteq [I:K]+[J:K]$. Since the reverse inclusion is always true, $[(I+J):K]=[I:K]+[J:K]$.

$(4) \rightarrow (5)$: Let $I$ and $J$ be finitely generated ideals of $S$. Then, $$S = [I+J:I+J] = [I:I+J]+[J:I+J] \subseteq [I:J]+[J:I] \subseteq S.$$

$(5) \rightarrow (6)$: (\cite[Proposition 4]{S}) It is clear that $[K:I]+[K:J] \subseteq [K: I\cap J]$. Let $s\in S$ such that $s(I\cap J) \subseteq K$. By hypothesis, $S=[I:J]+[J:I]$. So there exist $t_1\in [I:J]$ and $t_2\in [J:I]$ such that $1=t_1+t_2$. This implies that $s=st_1+st_2$. Let $x\in I$. Then $t_2x\in J$. Therefore, $t_2x\in I\cap J$. Since $s(I\cap J)\subseteq K$, $st_2x\in K$. Thus $st_2\in [K:I]$. Now let $y\in J$. Then $t_1y\in  I$. Therefore, $t_1y\in I\cap J$. Thus $st_1\in [K:J]$. So finally we have $s\in [K:I]+[K:J]$. Therefore, $[K:I\cap J] \subseteq [K:I]+[K:J]$.

$(6) \rightarrow (1)$: The proof is just a mimic of the proof of \cite[Theorem 6.6]{LarsenMcCarthy1971} and therefore, it is omitted.
\end{proof}

\end{theorem}

We end this section by characterizing Pr\"{u}fer semirings in terms of valuation semidomains. Note that valuation semirings have been introduced and investigated in \cite{Nasehpour2018}. Let us recall that a semiring is called to be a B\'{e}zout semiring if each of its finitely generated ideal is principal.

\begin{proposition}
	
	\label{bezoutisvaluation}
	
	A local semidomain is a valuation semidomain if and only if it is a B\'{e}zout semidomain.
	
	\begin{proof}
		$(\rightarrow):$ Straightforward.
		
		$(\leftarrow):$ Let $S$ be a local semidomain. Take $x,y \in S$ such that both of them are nonzero. Assume that $(x,y) = (d)$ for some nonzero $d\in S$. Define $x^{\prime} = x/d$ and $y^{\prime} = y/d$. It is clear that there are $a,b \in S$ such that $ax^{\prime} + by^{\prime} =1$. Since $S$ is local, one of $ax^{\prime}$ and $by^{\prime}$ must be unit, say $ax^{\prime}$. So $x^{\prime}$ is also unit and therefore, $(y^{\prime}) \subseteq S = (x^{\prime})$. Now multiplying the both sides of the inclusion by $d$ gives us the result $(y) \subseteq (x)$ and by Theorem 2.4 in \cite{Nasehpour2018}, the proof is complete.
	\end{proof}
	
\end{proposition}

Now we get the following nice result:

\begin{theorem}

\label{valuationlocalization}

For a semidomain $S$, the following statements are equivalent:

\begin{enumerate}

\item $S$ is Pr\"{u}fer.

\item For any prime ideal $\mathfrak{p}$, $S_\mathfrak{p}$ is a valuation semidomain.

\item For any maximal ideal $\mathfrak{m}$, $S_\mathfrak{m}$ is a valuation semidomain.

\end{enumerate}

\begin{proof}

$(1) \rightarrow (2):$

Let $J$ be a finitely generated nonzero ideal in $S_\mathfrak{p}$, generated by $s_1 / u_1, \ldots, s_n / u_n$, where $s_i \in S$ and $u_i \in S-\mathfrak{p}$. It is clear that $J = I_\mathfrak{p}$, where $I = (s_1, \ldots, s_n)$. By hypothesis, $I$ is invertible. So by Theorem \ref{invertiblelocalization2}, $J$ is principal. This means that $S_\mathfrak{p}$ is a B\'{e}zout semidomain and since it is local, by Proposition \ref{bezoutisvaluation}, $S_\mathfrak{p}$ is a valuation semidomain.

$(2) \rightarrow (3):$ Trivial.

$(3) \rightarrow (1):$ Let $I$ be a nonzero finitely generated ideal of $S$. Then for any maximal ideal $\mathfrak{m}$ of $S$, $I_\mathfrak{m}$ is a nonzero principal ideal of $S_\mathfrak{m}$ and by Theorem \ref{invertiblelocalization2}, $I$ is invertible. So we have proved that the semiring $S$ is Pr\"{u}fer and the proof is complete.
\end{proof}

\end{theorem}

Now we pass to the next section that is on Gaussian semirings.

\section{Gaussian Semirings}

In this section, we discuss Gaussian semirings. For doing so, we need to recall the concept of the content of a polynomial in semirings. Let us recall that for a polynomial $f \in S[X]$, the content of $f$, denoted by $c(f)$, is defined to be the finitely generated ideal of $S$ generated by the coefficients of $f$. In \cite[Theorem 3]{Nasehpour2016}, the semiring version of the Dedekind-Mertens lemma (Cf. \cite[p. 24]{Prufer1932} and \cite{ArnoldGilmer1970}) has been proved. We state that in the following only for the convenience of the reader:

\begin{theorem}[Dedekind-Mertens Lemma for Semirings]

\label{dedekindmertens2}

Let $S$ be a semiring. Then the following statements are equivalent:

\begin{enumerate}

\item The semiring $S$ is subtractive, i.e. each ideal of $S$ is subtractive,
\item If $f, g \in S[X]$ and $\deg(g) = m$, then $c(f)^{m+1} c(g) = c(f)^m c(fg)$.

\end{enumerate}

\end{theorem}

Now, we recall the definition of Gaussian semirings:

\begin{definition}
A semiring $S$ is said to be Gaussian if $c(fg)=c(f)c(g)$ for all polynomials $f,g \in S[X]$ \cite[Definition 7]{Nasehpour2016}.
\end{definition}

Note that this is the semiring version of the concept of Gaussian rings defined in \cite{Tsang1965}. For more on Gaussian rings, one may refer to \cite{BazzoniGlaz2007} also.

\begin{remark}

There is a point for the notion of Gaussian semirings that we need to clarify here. An Abelian semigroup $G$ with identity, satisfying the cancellation law, is called a Gaussian semigroup if each of its elements $g$, which is not a unit, can be factorized into the product of irreducible
elements, where any two such factorizations of the element $g$ are associated with each other \cite[\S 8 p. 71]{Kurosh1965}. In the papers \cite{DalePitts1978} and \cite{HebischWeinert1987} on Euclidean semirings, a semiring $S$ is called to be Gaussian if its semigroup of nonzero elements is Gaussian, which is another notion comparing to ours.

\end{remark}

Finally, we emphasize that by Theorem \ref{dedekindmertens2}, each ideal of a Gaussian semiring needs to be subtractive. Such semirings are called subtractive. Note that the boolean semiring $\mathbb B = \{0,1\}$ is a subtractive semiring, but the semiring $\mathbb N_0$ is not, since its ideal $\mathbb N_0 - \{1\}$ is not subtractive. As a matter of fact, all subtractive ideals of the semiring $\mathbb N_0$ are of the form $k \mathbb N_0$ for some $k\in \mathbb N_0$ \cite[Proposition 6]{Noronha1978}.

With this background, it is now easy to see that if every nonzero finitely generated ideal of a subtractive semiring $S$, is invertible, then $S$ is Gaussian. Also note that an important theorem in commutative ring theory, known as Gilmer-Tsang Theorem (cf. \cite{Gilmer1967} and \cite{Tsang1965}), states that $D$ is a Pr\"{u}fer domain if and only if $D$ is a Gaussian domain. The question may arise if a semiring version for Gilmer-Tsang Theorem can be proved. This is what we are going to do in the rest of the paper. First we prove the following interesting theorem:

\begin{theorem}

\label{Gaussian1}

 Let $S$ be a semiring. Then the following statements are equivalent:

 \begin{enumerate}

 \item $S$ is a Gaussian semidomain and $ab\in (a^2,b^2)$ for all $a,b\in S$,
 \item $S$ is a subtractive and Pr\"{u}fer semiring.

 \end{enumerate}

 \begin{proof}

 $(1) \rightarrow (2)$: Since $ab\in (a^2,b^2)$, there exists $r,s\in S$ such that $ab=ra^2+sb^2$. Now define $f,g\in S[X]$ by $f=a+bX$ and $g=sb+raX$. It is easy to check that $fg=sab+abX+rabX^2$. Since $S$ is Gaussian, $S$ is subtractive by Theorem \ref{dedekindmertens2}, and we have $c(fg)=c(f)c(g)$, i.e. $(ab)=(a,b)(sb,ra)$. But $(ab)=(a)(b)$ is invertible and therefore, $(a,b)$ is also invertible and by Theorem \ref{prufer}, $S$ is a Pr\"{u}fer semiring.

 $ (2) \rightarrow (1)$: Since $S$ is a subtractive and Pr\"{u}fer semiring, by Theorem \ref{dedekindmertens2}, $S$ is a Gaussian semiring. On the other hand, one can verify that $(ab)(a,b) \subseteq (a^2,b^2)(a,b)$ for any $a,b \in S$. If $a=b=0$, then there is nothing to be proved. Otherwise, since $(a,b)$ is an invertible ideal of $S$, we have $ab\in (a^2,b^2)$  and this completes the proof.
 \end{proof}

\end{theorem}

\begin{theorem}[Gilmer-Tsang Theorem for Semirings]

\label{Gilmer-Tsang}

Let $S$ be a subtractive semidomain such that $ab\in (a^2,b^2)$ for all $a,b \in S$. Then the following statements are equivalent:

 \begin{enumerate}

 \item $S$ is a Pr\"{u}fer semiring,
 \item Each nonzero finitely generated ideal of $S$ is cancellation,
 \item $[IJ \colon I]=J$ for all finitely generated ideals $I$ and $J$ of $S$,
 \item $S$ is a Gaussian semiring.
 \end{enumerate}

 \begin{proof}
 Obviously $(1) \rightarrow (2)$ and $(2) \rightarrow (3)$ hold by Proposition \ref{can1} and Proposition \ref{can2}, respectively.

 $(3) \rightarrow (4)$: Let $f,g \in S[X]$. By Theorem \ref{dedekindmertens2}, we have $c(f)c(g) c(f)^{m}=  c(fg) c(f)^m$. So $[c(f)c(g) c(f)^{m} \colon c(f)^{m}] =  [c(fg) c(f)^m \colon c(f)^{m}]$. This means that $c(f)c(g) =  c(fg)$ and $S$ is Gaussian.

 Finally, the implication $(4) \rightarrow (1)$ holds by Theorem \ref{Gaussian1} and this finishes the proof.
 \end{proof}

\end{theorem}

\begin{remark}

In \cite[Theorem 9]{Nasehpour2016}, it has been proved that every bounded distributive lattice is a Gaussian semiring. Also, note that if $L$ is a bounded distributive lattice with more than two elements, it is neither a ring nor a semidomain, since if it is a ring then the idempotency of addition causes $L=\{0\}$ and if it is a semidomain, the idempotency of multiplication causes $L=\mathbb B =\{0,1\}$. With the help of the following theorem, we give a plenty of examples of proper Gaussian and Pr\"{u}fer semirings. Let us recall that if $S$ is a semiring, then by $\FId(S)$, we mean the semiring of finitely generated ideals of $S$.

\end{remark}

\begin{theorem}

\label{notring}

Let $S$ be a Pr\"{u}fer semiring. Then the following statements hold for the semiring $\FId(S)$:

\begin{enumerate}

\item $\FId(S)$ is a Gaussian semiring.

\item $\FId(S)$ is a subtractive semiring.

\item $\FId(S)$ is an additively idempotent semidomain and for all finitely generated ideals $I$ and $J$ of $S$, we have $IJ \in (I^2,J^2)$.

\item $\FId(S)$ is a Pr\"{u}fer semiring.

\end{enumerate}

\begin{proof}

(1): Let $I,J \in \FId(S)$. Since $S$ is a Pr\"{u}fer semiring and $I \subseteq I+J$, by Theorem \ref{frac2}, there exists an ideal $K$ of $S$ such that $I=K(I+J)$. On the other hand, since $I$ is invertible, $K$ is also invertible. This means that $K$ is finitely generated and therefore, $K \in \FId(S)$ and $I \in (I+J)$. Similarly, it can be proved that $J \in (I+J)$. So, we have $(I,J) = (I+J)$ and by \cite[Theorem 8]{Nasehpour2016}, $\FId(S)$ is a Gaussian semiring.

(2): By Theorem \ref{dedekindmertens2}, every Gaussian semiring is subtractive. But by (1), $\FId(S)$ is a Gaussian semiring. Therefore, $\FId(S)$ is a subtractive semiring.

(3): Obviously $\FId(S)$ is additively-idempotent and since $S$ is a Pr\"{u}fer semiring, $\FId(S)$ is an additively idempotent semidomain. By Theorem \cite[Proposition 4.43]{Golan1999(b)}, we have $(I+J)^2 = I^2+J^2$ and so $(I+J)^2 \in (I^2,J^2)$. But $(I+J)^2 = I^2+J^2+IJ$ and by (2), $\FId(S)$ is subtractive. So, $IJ \in (I^2,J^2)$, for all $I,J \in \FId(S)$.

(4): Since $\FId(S)$ is a Gaussian semidomain such that $IJ\in (I^2,J^2)$ for all $I,J\in \FId(S)$, by Theorem \ref{Gaussian1}, $\FId(S)$ is a Pr\"{u}fer semiring and this is what we wanted to prove.
\end{proof}

\end{theorem}

\begin{corollary}

\label{notring2}

If $D$ is a Pr\"{u}fer domain, then $\FId(D)$ is a Gaussian and Pr\"{u}fer semiring.
\end{corollary}

\section*{Acknowledgments} The authors are very grateful to the anonymous referee for her/his useful pieces of advice, which helped them to improve the paper. The first named author is supported by the Faculty of Mathematics at the K. N. Toosi University of Technology. The second named author is supported by Department of Engineering Science at University of Tehran and Department of Engineering Science at Golpayegan University of Technology and his special thanks go to the both departments for providing all necessary facilities available to him for successfully conducting this research.

\bibliographystyle{plain}

\end{document}